 \documentclass{article}
 \usepackage[english]{babel}
 \usepackage[pdftex]{graphicx}
 \usepackage{wrapfig,multirow,fullpage,tensor,subfigure}
   \usepackage{fancybox,enumerate,floatrow,wrapfig,multirow}
  \usepackage{amsfonts,amsmath,amsthm,array,amssymb}
\usepackage{hyperref,color,adjustbox,lipsum,mathrsfs}
 \newtheorem{theorem}{Theorem}
 \include{latexsymb}
 \include{pxfonts}
 \include{displaymath}
 \include{tipa}

\begin{document}

 \begin{center}
 {\Large \textbf{On compatibility/incompatibility of two discrete probability distributions in the presence of incomplete specification}}
 \end{center}
 \bigskip
 \begin{center}
 {\bf \large Indranil Ghosh}\\ {\large University of North Carolina, Wilmington, North Carolina}\\
 e-mail: {\it  ghoshi@uncw.edu}\\
 \mbox{}\\
 \mbox{}\\
  {\bf \large N. Balakrishnan}\\ {\large McMaster University, Hamilton, Ontario, Canada}\\
 e-mail: {\it bala@mcmaster.ca}
 \end{center}

\bigskip
\begin{abstract}
Conditional specification of distributions is a developing area with many applications.
In the finite discrete case, a variety of compatible conditions can be derived.
In this paper, we propose an alternative approach to study the compatibility
of two conditional probability distributions under the finite discrete set up.
A technique based on rank--based criterion is shown to be particularly convenient
for identifying compatible distributions corresponding to complete conditional specification,  including the case with zeros.
The proposed methods are finally illustrated with several examples.
\end{abstract}
\bigskip


\section{Introduction}
Specification of joint distributions by means of conditional densities has received considerable attention in the literature in the last decade or so. Possible applications may be found in the area of model building  and in the elicitation and construction of multiparameter prior distributions in Bayesian scenarios. For example, suppose  $\underline{X}=(X_{1},X_{2},\cdots, X_{k})$ is a $k$-dimensional random vector taking on values in the finite range set $\underline{\mathcal{X}}_{1}\times \underline{\mathcal{X}}_{2} \times \cdots \times \underline{\mathcal{X}}_{k},$ where $\underline{\mathcal{X}}_{i}$ denotes the possible values of $X_{i},i=1,2, \cdots, k.$ Efforts to ascertain an appropriate distribution for $\underline{X}$ frequently involve acceptance or rejection of a series of bets about the stochastic behavior of $\underline{X}.$ Let us consider that in this situation  we are facing a question of whether or not to accept with odds 4 to 1 a bet that $X_{1}$ is equal to $1.$ Then, if we accept the bet  it puts a bound on the probability that $X=1.$

The basic problem is most easily visualized in the finite discrete case. Several different approaches exist in the literature with regard to the problem of determination of the possible compatibility of two families of conditional distributions (Arnold and Press, 1989; Arnold and Gokhale, 1994; Cacoullos and Papageorgiou, 1983; Wesolowski, 1995). In addition, the problem of determining most nearly compatible distributions, in the absence of compatibility, has been addressed [Arnold and Gokhale 1998; Arnold, Castillo and Sarabia (1999, 2001)]. Here we focus, on the finite discrete case, and take a closer look at the compatibility problem viewing it as a problem involving linear equations in restricted domains. The issue of near compatibility is also discussed using the concept of $\varepsilon$-compatibility [see Arnold et al. (1999), Ghosh and Balakrishnan (2013), and the references therein].  Furthermore, we also focus our attention on situations when we have incomplete (or partial) information on (either or both)  the two conditional probability matrices $A$ and $B,$ under the compatible set-up.

In particular,  we transform the problem of compatibility to a  linear programming problem and  derive conditions for compatibility based on the rank of a matrix $D,$ whose elements are functions of the two conditional probability matrices  $A$ and  $B.$ It is  found that the problem of compatibility, with our approach, is reduced to a large extent to a set of $IJ$ equations in $(I-1)$ unknowns with non-negativity constraints, where $I$ and $J$ are the dimensions of the matrices $A$ and $B.$ However,  we  mainly focus here on cases in which $I=2, 3, 4$ and $J=2, 3, 4.$  The rest of this paper  is organized as follows. In Section 2, we discuss  the concept of  compatibility for any two conditional probability matrices $A$ and $B$ in  the discrete set up.  In Section 3, we  discuss an alternative approach to  compatibility for the $(2\times 2)$ and $(3\times 3)$ cases and introduce the idea of rank--based criterion based on the rank of the matrix $D.$ In Section 4, we provide a discussion on  the problem of compatibility and/or minimal incompatibility when we have incomplete specification of matrices $A$ or $B,$ or both. In Section 5, some concluding remarks are provided.

\section{ Compatibility}
Let $A$ and $B$  be two  $(I\times J)$ matrices with non-negative elements such that $\sum_{i=1}^{I}a_{ij}=1$, $\forall j=1,\ldots,J$ and $\sum_{j=1}^{J} b_{ij} = 1$, $\forall i =1,2,\ldots,I$. Without loss of generality, it can be assumed that  $I\leq J.$
 Matrices  $A$ and  $B$ are said to form a compatible conditional specification for the distribution of ($X,Y$) if there exists some $(I\times J)$ matrix $P$ with non-negative entries  $p_{ij}$   and with $\sum_{i=1}^{I}\sum_{j=1}^{J}p_{ij}=1$ such that, for every $(i,j),$ $$a_{ij}=\frac{p_{ij}}{p_{.j}}$$ and $$b_{ij}=\frac{p_{ij}}{p_{i.}},$$ where $p_{i.}=\sum_{j=1}^{J}p_{ij}$ and $p_{i.}=\sum_{i=1}^{I}p_{ij}.$
 If such a matrix $P$ exists, then, if we assume that
 $$p_{ij}=P(X=x_{i}, Y=y_{j}),$$ $i=1,2,\cdots, I,$ $j=1,2,\cdots, J,$
 we will have
 $$a_{ij}=P(X=x_{i}|Y=y_{j}),$$ $i=1,2,\cdots, I,$ $j=1,2,\cdots, J,$
 and
 $$b_{ij}=P(Y=y_{j}|X=x_{i}),$$ $i=1,2,\cdots, I,$ $j=1,2,\cdots, J.$
Equivalently,  $A$ and  $B$ are compatible if there exist stochastic vectors  $\underline{\tau}=(\tau_{1},\tau_{2},\cdots,\tau_{J})$  and  $\underline{\eta}=(\eta_{1},\eta_{2},\cdots,\eta_{I})$  such that
\begin{equation*}
a_{ij}\tau_{j}=b_{ij}\eta_{i}
\end{equation*}
for every ($i,j$).
In the case of compatibility, $\underline{\eta}$ and $\underline{\tau}$ can be readily interpreted as the resulting  marginal distributions of $X$ and $Y,$ respectively.
 For any probability vector $\eta =({\eta_1,\eta_2,\ldots, \eta_I}),$  $p_{ij}=b_{ij}\eta_{i}$ is a probability distribution on the $IJ$ cells. So, the  conditional probability matrix, denoted by $A,$ and its elements $(a_{ij})$ will be given by
\begin{equation}
 a_{ij}=\frac{p_{ij}}{\displaystyle\sum_{s=1}^{I}p_{sj}}=\frac{b_{ij}\eta_i}{\displaystyle\sum_{s=1}^{I}b_{sj}\eta_s},
 \end{equation}
 for every ($i,j$).
If $A$ and $B$ are compatible, then
 $$a_{ij}\displaystyle\sum_{s=1}^{I}b_{sj} \eta_{s}=b_ {ij}\eta^{}_{i}.$$
 We then have  $$\tau_{j}=\displaystyle \sum_{s=1}^{I}b_{ij}\eta_{s},  \forall j=1,\ldots, J.$$
In this case, the  expressions given in   (1) can be rewritten as
\begin{equation*}
      a_{ij}\displaystyle\sum_{s=1}^{I}b_{sj}\eta^{}_{s}-b_{ij}\eta^{}_{i}=0,
  \end{equation*}

\noindent which in matrix notation the above can be written as
 \begin{equation}
 D\underline{\eta}=0,
 \end{equation}
where $D$ is a matrix of dimension $IJ\times I,$ and
the above equation is a system of  $IJ$  equations in $I-1$ unknowns  $\eta_{i},$ in view of the restriction $\sum_{i=1}^{I}\eta_{i}=1$. Through well-known matrix operations (such as left-multiplication by non-singular matrices), its rows can be reduced to at most  $I$ rows with non-zero elements (the so called $``\text{Row Echelon form}"$). Now, let this reduced system be denoted by  $D_{r}y = 0,$  where $y =(y_{1},y_{2},\ldots, y_{I-1})'$. Matrices $A$ and $B$ are compatible if the system $D_{r}y = 0$ has a solution $y$ of non-negative elements with at least one positive element. If such a $y^{*}$ exists, it can be scaled to arrive at a probability vector $\eta^{*}$. However, $A$ and $B$ are not compatible if the only solution with non-negative elements of  $D_{r}y=0$ is the null vector. In order to examine whether or not such a solution $y^{*}$ of  $D_{r}y = 0$  exists (especially when $I - 1$ is large), the methodology of linear programming may be used.  Specifically, consider the problem of
maximizing the objective function $\sum_{i}y_{i},$  subject to (a) the non-negativity constraints $\sum_{i}y_{i}\geq 0,$  (b) the equality constraints  $D_{r}y=0,$ and (c) the constraint $\sum_{i}y_{i}\leq 1$.  If the maximum of the objective function is positive, then the corresponding optimizing vector is $y^{*},$ which can be scaled into a probability vector $\eta^{*}$.  If the maximum is 0, then $A$ and $B$ are not compatible.

\subsection{ Compatibility of two matrices $A$ and $B$}
We know that if the matrices $A$ and $B$ are compatible, then $a_{ij}p_{\cdot j}=b_{ij}p_{i \cdot}$
for every $i = 1, 2, \ldots, I, j = 1, 2, \ldots, J$
(see Arnold et al. (1999)).
Equivalently, we can write
\begin{eqnarray*}
\displaystyle
a_{ij}\displaystyle\sum_{s=1}^{I}p_{sj}-b_{ij}\displaystyle\sum_{k=1}^{J}p_{ik}=0
\end{eqnarray*}
for every $i = 1, 2, \ldots, I, j = 1, 2, \ldots, J$,
which again can be written as
\begin{eqnarray*}
\displaystyle
a_{ij}\left[p_{1j}+p_{2j}+\cdots+p_{ij}+\cdots+p_{Ij}\right] -
b_{ij}\left[p_{i1}+p_{i2}+\cdots+p_{ij}+\cdots+p_{iJ}\right]=0
\end{eqnarray*}
for every $i = 1, 2, \ldots, I, j = 1, 2, \ldots, J$.
In matrix notation, the above system of linear equations can be written as
\begin{eqnarray}
\displaystyle
C\underline{p}=0,
\label{111}
\end{eqnarray}
where $C$ contains elements calculated from those of $A$ and $B$ and is a matrix of dimension $IJ\times IJ$ and
$\underline{p}^{(IJ\times 1)} = \left(p_{11},p_{12},\ldots,p_{IJ}\right)^{T}$.

\bigskip

\begin{theorem}
The solution space,  $\Omega$, for the system of equations in  (3)
is $(I-M)\underline{z}$,  where $M$ is an idempotent matrix and $\underline{z}^{(IJ\times 1)}$
is any arbitrary vector of dimension $IJ\times 1$.
\end{theorem}

\smallskip

\noindent \textit{Proof.} See Ghosh and Nadarajah (2017).

\section{An alternative approach to compatibility}
Questions  of compatibility of conditional and marginal specifications of distributions are of fundamental
importance in modeling. The earliest work in this regard is Patil (1965).
He considered the discrete case under a mild regularity condition and showed that the conditional distributions
of $X$ given $Y$ and of $Y$ given $X$ will uniquely determine the joint distribution  of $(X, Y)$.
There are several versions of necessary and sufficient conditions for compatibility given
by Arnold and Press (1989) and Arnold et al. (2002, 2004).
In some situations, the condition of  Arnold et al. (2004)  for checking compatibility
was found to be difficult and less effective.
This is the reason for us to propose here  a relatively easy and simple procedure to check the condition for compatibility.
The new method, which requires only some elementary type operation of matrices (``Row Echelon form''),
 provides a much simpler and an effective approach.

\noindent When the given conditional distributions are compatible,
it is natural to ask whether the associated joint distribution is unique.
This issue has been addressed in the literature by Amemiya (1975),
Gourieroux and Montfort (1979), Nerlove and Press (1986), and Arnold and Press (1989).
Arnold and Press (1989) pointed out that the condition for uniqueness is generally
difficult to check.
In this paper, through the structure of the reduced $D$ matrix,
we provide a simple criteria for checking uniqueness as well.

\noindent Here, we discuss the compatibility of two conditional matrices $A$ and $B$ along with the uniqueness
and the existence of a joint probability $P$ based on the rank of a matrix $D$.
A key  feature is  in the fact that it can be applied to situations
wherein matrices $A$ and $B$ have some zeros appearing in the same position.
In situations like this, the cross product criterion can not be applied to check compatibility.

\begin{theorem}
Any two given conditional probability matrices $A$ and $B$ of dimension $(I\times J)$
are compatible if $\text{\rm rank} \left(D^{(IJ\times I)} \right) \leq I-1$ with equality when there
exists a unique solution for the unknown $\eta_{i},$ for every $i.$
\end{theorem}

\smallskip

\noindent
{\bf Proof:}
Note that $\text{\rm rank}\left(D^{(IJ\times I)}\right) \leq \text{\rm min}(IJ,I)=I$.
Now when $D$ has full rank, i.e., $\text{\rm rank}(D)=I$, the only
solution to $D\underline{\eta}=0$ is the null vector (trivial solution).
Thus, matrices $A$ and $B$ are incompatible.
Next, if we have  $\text{\rm rank} \left( D^{(IJ\times I)} \right) \leq I-1$, the number
of equations $(IJ)$ is greater than the number of unknowns $(I-1),$ and so we must have a non-trivial solution.
If the non-trivial solution is positive then it can be appropriately scaled to arrive
at a probability vector $\underline{\eta}^{*}$.
Hence, the two matrices $A$ and $B$ are compatible.
However, in this case, the system of equations is not homogeneous and we have at most $(I-1)$ solutions.
When  $\text{\rm rank}(D)=I-1,$ we have $(I-1)$ unknowns
subject to the linear constraint $\displaystyle\sum_{i=1}^{I}\eta_{i}=1$.
The $(I-1)$ equations
(excluding the redundant equations from the total set of $IJ$ equations) and the system
of linear equations is homogeneous so that there exists a unique solution.
This completes the proof.
$\square$

\smallskip
\noindent This theorem is useful in situations when the two conditional matrices $A$ and $B$ have
zeros as elements appearing in the same position in which case we can not guarantee the existence of a compatible
matrix $P$ by the cross--product ratio criterion.
Next, we discuss  the  compatibility for  $(3\times 3),$  and $(4\times 4)$  cases with some examples.

\subsection{Proof of Rank($D$)=2 when $A$ and $B$ are compatible in a  $(3\times 3)$ case}
First of all, the form of the $D$-matrix in a $(3\times 3)$ case is given by
\begin{center}
\begin{equation*}
D=\left(
\begin{array}{ccc}
 b_{11}(a_{11}-1)&  a_{11}b_{21}& a_{11}b_{31}\\
  b_{12}(a_{12}-1)& a_{12}b_{22}& a_{12}b_{32}\\
  b_{13}(a_{13}-1)& a_{13}b_{23}& a_{13}b_{33}\\
  a_{21}b_{11} & b_{21}(a_{21}-1) & a_{21}b_{31} \\
   a_{22}b_{12}& b_{22}(a_{22}-1)& a_{22}b_{32}\\
  a_{23}b_{13}& b_{23}(a_{23}-1) & a_{23}b_{33}\\
   a_{31}b_{11} & a_{31}b_{21} & b_{31}(a_{31}-1)\\
    a_{32}b_{12} & a_{32}b_{22} & b_{32}(a_{32}-1)\\
  a_{33}b_{13} &a_{33}b_{23}& b_{33}(a_{33}-1) \\
\end{array}
\right).
\end{equation*}
\end{center}

Note that if the matrices $A$ and $B$ are compatible, then all possible cross product ratio($A$)=cross product ratio($B$).
First of all, we apply the following elementary row operations:
\begin{itemize}
\item new(row1)=row 1+row 4+row 7
\item new(row2)=row 2+row 5+row 8
\item new(row3)=row 3+row 6+row 9,
\end{itemize}
so that matrix $D$ reduces to
\begin{center}
\begin{equation*}
D=\left(
\begin{array}{ccccccccc}
 0&0&0\\
  0&0&0\\
  0&0&0\\
  a_{21}b_{11} &b_{21}(a_{21}-1) &a_{21}b_{31} \\
   a_{22}b_{12}&b_{22}(a_{22}-1)&a_{22}b_{32}\\
  a_{23}b_{13}&b_{23}(a_{23}-1) &a_{23}b_{33}\\
   a_{31}b_{11} &a_{31}b_{21} &b_{31}(a_{31}-1)\\
    a_{32}b_{12} &a_{32}b_{22} & b_{32}(a_{32}-1)\\
  a_{33}b_{13} &a_{33}b_{23}&b_{33}(a_{33}-1) \\
\end{array}
\right).
\end{equation*}
\end{center}

Again, we perform the following elementary row and column operations:
\begin{itemize}
 \item new(row5)=$\frac{\text{row} 5}{a_{22}}$
\item new(row6)=$\frac{\text{row} 4}{a_{23}}$
\item new(row8)=$\frac{\text{row} 8}{a_{32}}$
\item new(row9)=$\frac{\text{row} 4}{a_{33}}$
\item new(row 4)=$\frac{\text{row} 4}{a_{21}}$+new(row 5)+new(row 6)
\item new(row7)=$\frac{\text{row} 7}{a_{31}}$+new(row 8)+new(row 9),
\end{itemize}

 \noindent so that  matrix  $D$ has the form
\begin{center}
\begin{equation*}
D=\left(
\begin{array}{ccc}
 0&0&0\\
  0&0&0\\
  0&0&0\\
  1 &1-(\frac{b_{21}}{a_{21}}+\frac{b_{22}}{a_{22}}+\frac{b_{23}}{a_{23}}) &1 \\
   b_{12}&b_{22}(1-\frac{1}{a_{22}})&b_{32}\\
  b_{13}&b_{23}(1-\frac{1}{a_{23}}) &b_{33}\\
   1 &1 &1-(\frac{b_{31}}{a_{31}}+\frac{b_{32}}{a_{32}}+\frac{b_{33}}{a_{33}})\\
    b_{12} &b_{22} & b_{32}(1-\frac{1}{a_{32}})\\
  b_{13} &b_{23}&b_{33}(1-\frac{1}{a_{33}}) \\
\end{array}
\right).
\end{equation*}
\end{center}
Now,  we consider the following row operations:
\begin{itemize}
\item new(row 8)=row 8+row 9
\item new(row 5)=2 row 5+row 6-new(row 8)
\item new(row 4)=row 4-row 7
\item new(row 6)=row 6-row 8,
\end{itemize}

with which matrix  $D$ reduces to
\begin{center}
\begin{equation*}
D=\left(
\begin{array}{ccc}
 0&0&0\\
  0&0&0\\
  0&0&0\\
  0 &-(\frac{b_{21}}{a_{21}}+\frac{b_{22}}{a_{22}}+\frac{b_{23}}{a_{23}}) &(\frac{b_{31}}{a_{31}}+\frac{b_{32}}{a_{32}}+\frac{b_{33}}{a_{33}}) \\
   0&-(\frac{b_{22}}{a_{22}}+\frac{b_{23}}{a_{23}})&(\frac{b_{32}}{a_{32}}+\frac{b_{33}}{a_{33}})\\
  0&-\frac{b_{23}}{a_{23}} &\frac{b_{33}}{a_{33}}\\
   1 &1 &1-(\frac{b_{31}}{a_{31}}+\frac{b_{32}}{a_{32}}+\frac{b_{33}}{a_{33}})\\
   - b_{11} &-b_{21} &-b_{31}-(\frac{b_{32}}{a_{32}}+\frac{b_{33}}{a_{33}}) \\
  b_{13} &b_{23}&b_{33}(1-\frac{1}{a_{33}})\\
\end{array}
\right).
\end{equation*}
\end{center}

Now, with new (row 4)=row 4-row 5,  new(row 5)=row 5-row 6,  the $D$ matrix reduces to
\begin{center}
\begin{equation*}
D=\left(
\begin{array}{ccccccccc}
 0&0&0\\
  0&0&0\\
  0&0&0\\
  0 &-\frac{b_{21}}{a_{21}} &\frac{b_{31}}{a_{31}} \\
   0&-\frac{b_{22}}{a_{22}}&\frac{b_{32}}{a_{32}}\\
  0&-\frac{b_{23}}{a_{23}} &\frac{b_{33}}{a_{33}}\\
   1 &1 &1-(\frac{b_{31}}{a_{31}}+\frac{b_{32}}{a_{32}}+\frac{b_{33}}{a_{33}})\\
   - b_{11} &-b_{21} &-b_{31}-(\frac{b_{32}}{a_{32}}+\frac{b_{33}}{a_{33}}) \\
  b_{13} &b_{23}&b_{33}(1-\frac{1}{a_{33}})\\
\end{array}
\right).
\end{equation*}
\end{center}

Let us consider the determinant of any sub-matrix of order $(2\times 2),$ say,

 $$B=\left(
\begin{array}{cc}
 -b_{11} &-b_{21}\\
 b_{13} &b_{23}.\\
\end{array}
\right)$$

The determinant for the matrix $B$ is given by
\begin{eqnarray}
\text{det}(B)&=&-b_{11}b_{23}+ b_{13}b_{21}\nonumber\\
&\neq& 0.
\end{eqnarray}

\noindent Thus, we have  $\text{rank(D)}=I-1=2.$
Hence, $A$ and $B$ are compatible iff rank($D$)=$I-1$. However, if $A$ and $B$ are not compatible, then rows of $A$ are not proportional
to the rows of $B,$ which implies that $\text{rank(D)}>2$.
This completes the proof.  $\square$

\subsection{Proof of Rank($D$)=3 when $A$ and $B$ are compatible in a  $(4\times 4)$ case}
As before, the form of the $D$-matrix in a $(4\times 4)$ case is given by
\begin{center}
\begin{equation*}
D=\left(
\begin{array}{cccc}
 b_{11}(a_{11}-1)&  a_{11}b_{21}& a_{11}b_{31}& a_{11}b_{41}\\
  b_{12}(a_{12}-1)& a_{12}b_{22}& a_{12}b_{32}& a_{12}b_{42}\\
  b_{13}(a_{13}-1)& a_{13}b_{23}& a_{13}b_{33}& a_{13}b_{43}\\
   b_{14}(a_{14}-1)& a_{14}b_{24}& a_{14}b_{34}& a_{14}b_{44}\\
  a_{21}b_{11} & b_{21}(a_{21}-1) & a_{21}b_{31}& a_{21}b_{41} \\
   a_{22}b_{12}& b_{22}(a_{22}-1)& a_{22}b_{32}& a_{22}b_{42}\\
  a_{23}b_{13}& b_{23}(a_{23}-1) & a_{23}b_{33}& a_{23}b_{43}\\
  a_{24}b_{14}& b_{24}(a_{24}-1) & a_{24}b_{34}& a_{24}b_{44}\\
   a_{31}b_{11} & a_{31}b_{21} & b_{31}(a_{31}-1)& a_{31}b_{41}\\
    a_{32}b_{12} & a_{32}b_{22} & b_{32}(a_{32}-1)& a_{32}b_{42}\\
  a_{33}b_{13} &a_{33}b_{23}& b_{33}(a_{33}-1)& a_{33}b_{43} \\
  a_{34}b_{14} &a_{34}b_{24}& b_{34}(a_{34}-1)& a_{34}b_{44} \\
  a_{41}b_{11} &a_{41}b_{21}&a_{41}b_{31}& b_{41}(a_{41}-1)  \\
  a_{42}b_{12} &a_{42}b_{22}&a_{42}b_{32}& b_{42}(a_{42}-1)  \\
  a_{43}b_{13} &a_{43}b_{23}&a_{43}b_{33}& b_{43}(a_{43}-1)  \\
  a_{44}b_{14} &a_{44}b_{24}&a_{43}b_{34}& b_{44}(a_{44}-1)\\
\end{array}
\right).
\end{equation*}
\end{center}

First, we consider the following elementary row operations:
\begin{itemize}
\item new(row1)=row 1+row 5+row 9+row 13
\item new(row2)=row 2+row 6+row 10+row 14
\item new(row3)=row 3+row 7+row 11+row 15
\item new(row3)=row 4+row 8+row 12+row 16,
\end{itemize}
so that matrix $D$ reduces to
\begin{center}
\begin{equation*}
D=\left(
\begin{array}{cccc}
 0&0&0&0\\
  0&0&0&0\\
  0&0&0&0\\
0&0&0&0\\
 a_{21}b_{11} & b_{21}(a_{21}-1) & a_{21}b_{31}& a_{21}b_{41} \\
   a_{22}b_{12}& b_{22}(a_{22}-1)& a_{22}b_{32}& a_{22}b_{42}\\
  a_{23}b_{13}& b_{23}(a_{23}-1) & a_{23}b_{33}& a_{23}b_{43}\\
  a_{24}b_{14}& b_{24}(a_{24}-1) & a_{24}b_{34}& a_{24}b_{44}\\
   a_{31}b_{11} & a_{31}b_{21} & b_{31}(a_{31}-1)& a_{31}b_{41}\\
    a_{32}b_{12} & a_{32}b_{22} & b_{32}(a_{32}-1)& a_{32}b_{42}\\
  a_{33}b_{13} &a_{33}b_{23}& b_{33}(a_{33}-1)& a_{33}b_{43} \\
  a_{34}b_{14} &a_{34}b_{24}& b_{34}(a_{34}-1)& a_{34}b_{44} \\
  a_{41}b_{11} &a_{41}b_{21}&a_{41}b_{31}& b_{41}(a_{41}-1)  \\
  a_{42}b_{12} &a_{42}b_{22}&a_{42}b_{32}& b_{42}(a_{42}-1)  \\
  a_{43}b_{13} &a_{43}b_{23}&a_{43}b_{33}& b_{43}(a_{43}-1)  \\
  a_{44}b_{14} &a_{44}b_{24}&a_{44}b_{34}& b_{44}(a_{44}-1) \\
\end{array}
\right).
\end{equation*}
\end{center}

\noindent Now, we consider the following elementary row and column operations:
\begin{itemize}
\item new(row 6)=$\frac{\text{row} 6}{a_{22}}$
\item new(row 7)=$\frac{\text{row} 7}{a_{23}}$
\item new(row 8)=$\frac{\text{row} 8}{a_{24}}$
\item new(row 10)=$\frac{\text{row} 10}{a_{32}}$
\item new(row 11)=$\frac{\text{row} 11}{a_{33}}$
\item new(row 12)=$\frac{\text{row} 12}{a_{34}}$
\item new(row 14)=$\frac{\text{row} 14}{a_{42}}$
\item new(row 15)=$\frac{\text{row} 15}{a_{43}}$
\item new(row 16)=$\frac{\text{row} 16}{a_{44}}$
\item new(row 5)=$\frac{\text{row}  5}{a_{21}}$+new(row 6)+new(row 7)+new(row 8)
\item new(row 9)=$\frac{\text{row}  9}{a_{31}}$+new(row 10)+new(row 11)+new(row 12)
\item new(row 13)=$\frac{\text{row} 13}{a_{41}}$+new(row 14)+new(row 15)+new(row 16),

\end{itemize}

 \noindent so that matrix $D$  has the form
\begin{center}
\begin{equation*}
D=\left(
\begin{array}{cccc}
 0&0&0&0\\
  0&0&0&0\\
  0&0&0&0\\
0&0&0&0\\
 1 & 1-\left(\frac{b_{21}}{a_{21}}+\frac{b_{22}}{a_{22}}+\frac{b_{23}}{a_{23}}+ \frac{b_{24}}{a_{24}}\right)& 1& 1 \\
   b_{12}& b_{22}\left(1-a^{-1}_{22}\right) & b_{32}& b_{42}\\
  b_{13}& b_{23}\left(1-a^{-1}_{23}\right) & b_{33}& b_{43}\\
  b_{14}& b_{24}\left(1-a^{-1}_{24}\right) & b_{34}& b_{44}\\
   1 & 1 & 1-\left(\frac{b_{31}}{a_{31}}+\frac{b_{32}}{a_{32}}+\frac{b_{33}}{a_{33}}+ \frac{b_{34}}{a_{34}}\right)& 1\\
    b_{12} & b_{22} & b_{32}\left(1-a^{-1}_{32}\right) & b_{42}\\
  b_{13} &b_{23}& b_{33}\left(1-a^{-1}_{33}\right) & b_{43} \\
  b_{14} &b_{24}& b_{34}\left(1-a^{-1}_{34}\right) & b_{44} \\
  1 &1&1&  1-\left(\frac{b_{41}}{a_{41}}+\frac{b_{42}}{a_{42}}+\frac{b_{43}}{a_{43}}+ \frac{b_{44}}{a_{44}}\right) \\
  b_{12} &b_{22}&b_{32}& b_{42}\left(1-a^{-1}_{42}\right)  \\
  b_{13} &b_{23}&b_{33}& b_{43}\left(1-a^{-1}_{43}\right)   \\
  b_{14} &b_{24}&b_{34}& b_{44}\left(1-a^{-1}_{44}\right)  \\
\end{array}
\right).
\end{equation*}
\end{center}
 We now consider the following row operations:
\begin{itemize}
\item new(row 5)=row 5-row 9
\item new(row 6)=row 6-row 10
\item new(row 7)=row 7-row 11
\item new(row 8)=row 8-row 12
\item new(row 9)=row 9-row 13
\item new(row 10)=row 10-row 14
\item new(row 11)=row 11-row 15
\item new(row 12)=row 12-row 16,
\end{itemize}

with which  matrix $D$  reduces to
\begin{center}
\begin{equation*}
D=\left(
\begin{array}{cccc}
 0&0&0&0\\
  0&0&0&0\\
  0&0&0&0\\
0&0&0&0\\
   0 & -\left(\frac{b_{21}}{a_{21}}+\frac{b_{22}}{a_{22}}+\frac{b_{23}}{a_{23}}+ \frac{b_{24}}{a_{24}}\right)& 1& 0 \\
   0& -\frac{b_{22}}{a_{22}} &  -\frac{b_{32}}{a_{32}}& 0\\
  0& -\frac{b_{23}}{a_{23}}&  -\frac{b_{33}}{a_{33}}& 0\\
  0& -\frac{b_{24}}{a_{24}} &  -\frac{b_{34}}{a_{34}}& 0\\
   0 & 0& -\left(\frac{b_{31}}{a_{31}}+\frac{b_{32}}{a_{32}}+\frac{b_{33}}{a_{33}}+ \frac{b_{34}}{a_{34}}\right)&- \left(\frac{b_{41}}{a_{41}}+\frac{b_{42}}{a_{42}}+\frac{b_{43}}{a_{43}}+ \frac{b_{44}}{a_{44}}\right)\\
    0& 0&  -\frac{b_{32}}{a_{32}}&  -\frac{b_{42}}{a_{42}}\\
  0 &0&  -\frac{b_{33}}{a_{33}} &  -\frac{b_{43}}{a_{43}} \\
  0&0&  -\frac{b_{34}}{a_{34}} &  -\frac{b_{44}}{a_{44}} \\
  1 &1&1&  1-\left(\frac{b_{41}}{a_{41}}+\frac{b_{42}}{a_{42}}+\frac{b_{43}}{a_{43}}+ \frac{b_{44}}{a_{44}}\right) \\
  b_{12} &b_{22}&b_{32}& b_{42}\left(1-a^{-1}_{42}\right)  \\
  b_{13} &b_{23}&b_{33}& b_{43}\left(1-a^{-1}_{43}\right)   \\
  b_{14} &b_{24}&b_{34}& b_{44}\left(1-a^{-1}_{44}\right)  \\
  \end{array}
\right).
\end{equation*}
\end{center}

Now, let us consider the determinant of any sub-matrix of order $(3\times 3),$ say,

 $$M=\left(
\begin{array}{ccc}
 b_{12} &b_{22}&b_{32}  \\
  b_{13} &b_{23}&b_{33} \\
  b_{14} &b_{24}&b_{34} \\
  \end{array}
\right).$$

The determinant of matrix $M$ is given by
\begin{eqnarray}
\text{det}(M)&=&b_{12} \left(b_{23} b_{34} -b_{24} b_{23} \right)-b_{22} \left(b_{13} b_{34} -b_{14} b_{33} \right)+b_{32} \left(b_{13} b_{24} -b_{14} b_{23} \right)
\nonumber\\
&\neq& 0.
\end{eqnarray}

\noindent Thus, we have  $\text{rank(D)}=I-1=3.$
Therefore, $A$ and $B$ are compatible if and only if  rank($D$)=$I-1$. If $A$ and $B$ are not compatible, then rows of $A$ are not proportional
to the rows of $B,$ which implies that $\text{rank(D)}>3$.
This completes the proof. $\square$

\section{Study of compatibility under incomplete specification on $A$ or $B,$ or both}
In this section we will consider the problem of compatibility of two conditional probability matrices  $A$ and $B$ under the discrete set-up when more than one element either in $A$ or in $B$ is unknown. In particular, we will discuss in detail the $(2\times 3)$ case and we will consider two different situations which in detail as follows:
\begin{itemize}
\item More than one element is unknown only in $A,$
\item More than one element is unknown in both $A$ and $B.$
\end{itemize}

Our objective here is to investigate what happens to the compatibility condition when we have above situations.
\subsection{Compatibility when only elements of $A$ are unknown}
\begin{enumerate}
\item Let us consider $I=2$ and $J=3,$ and we assume that only two elements of $A$ are  unknown while all the elements of $B$ are known. We denote the $(i,j)-th$ unknown element of $A$ by $\alpha_{ij}.$ Suppose
$$A=\left(
\begin{array}{ccc}
 a_{11}&\alpha_{12}&a_{13}\\
  a_{21}&\alpha_{22}&a_{23}\\
    \end{array}
\right)$$
  and

$$B=\left(
\begin{array}{ccc}
 b_{11}&b_{12}&b_{13}\\
  b_{21}&b_{22}&b_{23}\\
    \end{array}
\right).$$

 \noindent  Here, we assume that all the elements in matrices $A$ and $B$ are strictly positive. Also, $A$ has elements such that column sums are equal to one and $B$ has elements such that the row sums are equal to one. So,  we have $$\alpha_{11}+\alpha_{22}=1.$$

 We know that the problem of  compatibility can be reduced to (in matrix notation as)
 $D\underline{\eta}=0,$  where $D$ has elements computed from the matrices $A$ and $B.$  Also, we note that if the two matrices $A$ and $B$ are compatible, then $C\underline{p}=0$ and vice versa. In this case, we have two  constraints $\alpha_{12}+\alpha_{22}=1$ and $\eta_{1}+\eta_{2}=1.$
So, the set of  equations, involving $\alpha_{12}$ and $\alpha_{22},$ that are sufficient to finding the unknown values (remaining equations will be redundant), from  (2), will be

\begin{gather}
b_{12}(\alpha_{12}-1)\eta_{1}+\alpha_{12}b_{22}\eta_{2}=0,\\
b_{22}(\alpha_{22}-1)\eta_{2}+\alpha_{22}b_{12}\eta_{1}=0,\\
b_{13}(a_{13}-1)\eta_{1}+a_{13}b_{23}\eta_{2}=0.
\end{gather}

Now, due to the  constraint, we get from (6) that
$$\eta_{1}=\frac{a_{13}b_{23}}{a_{13}b_{23}+b_{13}(1-a_{13})}.$$

Again, by substituting the value of $\eta_{1}$ in (6) and using the  constraint that  $\alpha_{12}+\alpha_{22}=1,$ we get the value of $\alpha_{22},$ to be
$$\alpha_{22}=\frac{b_{22}(1-\eta_{1})}{b_{22}(1-\eta_{1})+b_{12}\eta_{1}}
=\frac{b_{22}b_{13}(1-a_{13})}{b_{22}b_{13}(1-a_{13})+b_{22}b_{13}a_{13}}.$$
\noindent Subsequently, the unknown value of $\alpha_{12}$ will be $\alpha_{12}=1-\alpha_{22}.$

\bigskip

\noindent {\bf Some Examples}
\begin{enumerate}

 \item  Suppose we have two matrices $A$ and $B$  as follows:
 $$A=\left(
\begin{array}{ccc}
 1/5&\alpha_{12}&3/4\\
  4/5&\alpha_{22}&1/4\\
    \end{array}
\right)$$
  and

$$B=\left(
\begin{array}{ccc}
1/6&2/6&3/6\\
  4/6&1/6&1/6\\
    \end{array}
\right).$$

Now, if we are given  that $A$ and $B$ are compatible, then the values of $\alpha_{12}$ and $\alpha_{12}$ will be given by
$$\alpha_{22}=\frac{b_{22}b_{13}(1-a_{13})}{b_{22}b_{13}(1-a_{13})+b_{22}b_{13}a_{13}}
=\frac{\frac{1}{6}\frac{3}{6}(1-\frac{3}{4})}{\frac{1}{6}\frac{3}{6}(1-\frac{3}{4})+\frac{2}{6}\frac{1}{6}\frac{3}{4}}
=\frac{1}{3}.$$
So, $\alpha_{12}=1-\frac{1}{3}=\frac{2}{3}.$
Note that these are the unique choices for the unknown elements in the matrix $A$ for which the above matrices are compatible.
\item Next, we consider the situation when $I=3$ and $J=3$ and, as before, denoting the unknown values of the matrix $A$ by $\alpha_{ij},$ in the $(i,j)-th$ position we have (with all elements of $B$ being known), where the matrices $A$ and $B$ are of the form

$$A=\left(
\begin{array}{ccc}
 a_{11}&\alpha_{12}&a_{13}\\
  a_{21}&\alpha_{22}&a_{23}\\
    a_{31}&\alpha_{32}&a_{33}\\
    \end{array}
\right)$$
  and

$$B=\left(
\begin{array}{ccc}
 b_{11}&b_{12}&b_{13}\\
  b_{21}&b_{22}&b_{23}\\
    b_{31}&b_{32}&b_{33}\\
    \end{array}
\right).$$

 The linear constraints in this case are as follows (considering the fact that the column sums of the matrix $A$ are each equal to one and $\eta_{i}, i=1,2,3,$ are the marginal probability vectors corresponding to $B$)
\begin{gather}
\alpha_{12}+\alpha_{22}+\alpha_{32}=1,\\
 \eta_{1}+\eta_{2}+\eta_{3}=1.
 \end{gather}
Then, according to the compatibility condition, we will have $D\underline{\eta}=0 $ if matrices $A$ and $B$ are compatible. However,
    the $D$ matrix in this case will be
    \begin{center}
    \begin{equation*}
\begin{bmatrix}
b_{11}(a_{11}-1)&  a_{11}b_{21}& a_{11}b_{31}\\
  b_{12}(\alpha_{12}-1)& \alpha_{12}b_{22}& \alpha_{12}b_{32}\\
  b_{13}(a_{13}-1)& a_{13}b_{23}& a_{13}b_{33}\\
  a_{21}b_{11} & b_{21}(a_{21}-1) & a_{21}b_{31} \\
  \alpha_{22} b_{12}& b_{22}(\alpha_{22}-1)& \alpha_{22}b_{32}\\
  a_{23}b_{13}& b_{23}(a_{23}-1) & a_{23}b_{33}\\
   a_{31}b_{11} & a_{31}b_{21} & b_{31}(a_{31}-1)\\
   \alpha_{32} b_{12} & \alpha_{32}b_{22} & b_{32}(\alpha_{32}-1)\\
  a_{33}b_{13} &a_{33}b_{23}& b_{33}(a_{33}-1) \\
\end{bmatrix}.
\end{equation*}
\end{center}
\noindent So, the set of linear equations to find the unknown $\eta_{i}'s$ as well as the unknown $\alpha_{ij}'s$ will be (from the above $D$ matrix) as follows:
\begin{gather}
b_{11}(a_{11}-1)\eta_{1}+ a_{11}b_{21}\eta_{2}+ a_{11}b_{31}\eta_{3}=0,\\
b_{13}(a_{13}-1)\eta_{1}+ a_{13}b_{23}\eta_{2}+ a_{13}b_{33}\eta_{3}=0,\\
b_{11}a_{21}\eta_{1}+ (a_{21}-1)b_{21}\eta_{2}+ a_{21}b_{31}\eta_{3}=0,\\
b_{12}(\alpha_{12}-1)\eta_{1}+ \alpha_{12}b_{22}\eta_{2}+ \alpha_{12}b_{32}\eta_{3}=0,\\
b_{12}\alpha_{22}\eta_{1}+ (\alpha_{22}-1)b_{22}\eta_{2}+ \alpha_{22}b_{32}\eta_{3}=0,\\
b_{12}\alpha_{32}\eta_{1}+ \alpha_{32}b_{22}\eta_{2}+ (\alpha_{32}-1)b_{32}\eta_{3}=0.
\end{gather}
Solving the above set of equations with the constraints,  we get the following expressions for the unknowns:
\begin{gather}
\eta_{1}=d_{22},\\
\eta_{2}=\frac{d_{11}}{d_{12}}=d_{13},\\
\eta_{3}=1-d_{22}-d_{13},\\
\alpha_{12}=\frac{b_{12}d_{22}}{b_{12}d_{22}+b_{22}d_{13}+b_{32}d_{23}},\\
\alpha_{22}=\frac{b_{12}d_{13}}{b_{12}d_{22}+b_{12}d_{13}+b_{32}d_{23}},\\
\alpha_{32}=\frac{b_{32}d_{23}}{b_{12}d_{22}+b_{22}d_{13}+b_{32}d_{23}},
\end{gather}
where
\begin{gather}
 d_{11}=a_{11}b_{31}[a_{13}b_{33}-b_{13}(a_{13}-1)]+a_{13}b_{33}[b_{11}(a_{11}-1)-a_{11}b_{31}],\\
 d_{12}=[(a_{11}b_{21}-a_{11}b_{31})(b_{13}(a_{13}-1)-a_{13}b_{33})]-[(a_{13}b_{23}-a_{13}b_{33})(b_{11}(a_{11}-1)-a_{11}b_{31})],\\
 d_{22}=\frac{d_{13}(a_{11}b_{31}-a_{11}b_{21})-a_{11}b_{31}}{b_{11}(a_{11}-1)-a_{11}b_{31}}.
 \end{gather}
 \noindent Consequently, on using (23), (24) in (18), one can get an expression of $d_{13}.$

 Next, let us consider, a situation for example, the following choices for the two matrices $A$ and $B:$
 $$A=\left(
\begin{array}{ccc}
 1/6&\alpha_{12}&1/4\\
  1/3&\alpha_{22}&7/16\\
    3/6&\alpha_{32}&5/16\\
    \end{array}
\right)$$
  and

$$B=\left(
\begin{array}{ccc}
 1/7&2/7&4/7\\
  2/5&2/5&1/5\\
    1/4&1/4&2/4\\
    \end{array}
\right).$$
\noindent Hence, using  Eqs. (23)-(25), we obtain
$d_{11}=-0.002976190,d_{12}=-0.0142857,d_{13}=\frac{d_{11}}{d_{12}}=0.208334, d_{22}=0.2916667.$
So, the unknown elements of the matrix $A$ will be
$$\alpha_{12}=0.2857165, \alpha_{22}=0.222248, \alpha_{32}=0.4920587.$$
\noindent Importantly, these are the unique choices for which the two given matrices $A$ and $B$ are compatible.
 \end{enumerate}

\subsection{Compatibility when some elements in both  $A$ and $B$ are unknown}
  Suppose  we have a situation  where in both in $A$ and  $B$ some  elements are unknown, and we define the unknown elements of the matrix $A$ by $\alpha_{ij}$ and unknown elements of the matrix $B$ by $\beta_{ij}.$
 First, let us  consider the situation when $I=2$ and $J=3$ with
  $$A=\left(
\begin{array}{ccc}
   a_{11}&\alpha_{12}&a_{11}\\
  a_{21}&\alpha_{22}&a_{23}\\
  \end{array}
\right)$$ and

  $$B=\left(
\begin{array}{ccc}
 b_{11}&\beta_{12}&\beta_{13}\\
  b_{21}&b_{22}&b_{23}\\
    \end{array}
\right).$$

\noindent  In this case, we have same  constraints on the unknown elements $\alpha_{ij}$ as before, and for  $\beta_{ij}$ we have the following  restrictions:
\begin{gather}
b_{11}+ \beta_{12}+\beta_{13}=1,\\
\alpha_{12}+\alpha_{22}=1,\\
\eta_{1}+\eta_{2}=1.
\end{gather}

  \noindent We will then have the following set of  equations (for those involving the unknowns and excluding the redundant equations):
\begin{gather}
b_{11}(a_{11}-1)\eta_{1}+a_{11}b_{21}\eta_{2}=0,\\
\beta_{12}(\alpha_{12}-1)\eta_{1}+\alpha_{12}b_{22}\eta_{2}=0,\\
b_{22}(\alpha_{22}-1)\eta_{2}+\alpha_{22}\beta_{12}\eta_{1}=0,\\
\beta_{13}(a_{13}-1)\eta_{1}+a_{13}b_{23}\eta_{2}=0.
\end{gather}
\noindent Again, by using the  constraints,  we get from (32) that,
\begin{equation}
\beta_{13}=\frac{(1-\eta_{1})a_{13}b_{23}}{(1-a_{13})\eta_{1}}.
\end{equation}

\noindent Also, from  (29), by using the  constraint in (28),  we get
\begin{equation*}
\eta_{1}=\frac{a_{11}b_{21}}{a_{11}b_{21}+b_{11}(1-a_{11})}.
\end{equation*}

\noindent Substituting the above expression of $\eta_{1}$ in (33), we get (after some algebraic simplification),

\begin{equation}
\beta_{13}=\frac{ b_{11}b_{23}a_{21}a_{13}}{a_{11}a_{23}b_{21}}
\end{equation}

\noindent Hence, the  value of $\beta_{12}$ becomes

\begin{equation}
\beta_{12}=1-b_{11}-\frac{ b_{11}b_{23}a_{21}a_{13}}{a_{11}a_{23}b_{21}}.
\end{equation}

\noindent Again, substituting in  (30), we get

\begin{eqnarray}
\alpha_{12}&=&\frac{b_{21}a_{11}\beta_{12}}{b_{21}a_{11}b_{22}+b_{21}a_{11}\beta_{12}}\nonumber\\
&=&\frac{(a_{11}a_{23}b_{21}- b_{11}b_{23}a_{21}a_{13})a_{11}b_{21}}
{(a_{11}a_{23}b_{21}- b_{11}b_{23}a_{21}a_{13})a_{11}b_{21}+b_{11}b_{22}a_{21}a_{11}a_{23}b_{21}}.
\end{eqnarray}

Due to the  constraint, we can now find the unknown value of $\alpha_{22}$  to be
$$\alpha_{22}=1-\alpha_{12}.$$

As an example,  As before let us consider matrices $A$ and $B$ as
$$A=\left(
\begin{array}{ccc}
 1/5&\alpha_{12}&1/2\\
  4/5&\alpha_{22}&1/2\\
    \end{array}
\right)$$
  and

$$B=\left(
\begin{array}{ccc}
1/6&\beta_{12}&\beta_{13}\\
  2/5&2/5&1/5\\
    \end{array}
\right).$$

Here, if we are given the information that $A$ and $B$ are compatible, then the choices of the unknown values of $\alpha_{ij}'s$ and $\beta_{ij}'s$ will be given by
\begin{equation*}
\beta_{13}=\frac{ b_{11}b_{23}a_{21}a_{13}}{a_{11}a_{23}b_{21}}
=\frac{1}{3},
\end{equation*}

\begin{equation*}
\beta_{12}=1-\frac{1}{6}-\frac{1}{3}=\frac{1}{2},
\end{equation*}

\begin{equation*}
\eta_{1}=\frac{\frac{1}{5}\frac{2}{5}}{\frac{1}{5}\frac{2}{5}+\frac{1}{6}\frac{4}{5}}
=\frac{3}{8},
\end{equation*}

so that
\begin{eqnarray*}
\alpha_{12}&=&\frac{\beta_{12}\eta_{1}}{\beta_{12}\eta_{1}+b_{22}(1-\eta_{1})}\nonumber\\
&=&\frac{3}{7},
\end{eqnarray*}

 and hence $$\alpha_{22}=1-\alpha_{12}=\frac{4}{7}.$$
Furthermore, we note that in this case also, these are the choices for the unknown values for which matrices $A$ and $B$ are compatible and that they are unique.

\subsection{Choices of the unknown values of $A$ under incompatibility}
  Let us define $d_{ij}=\frac{a_{ij}}{b_{ij}},$ provided $b_{ij}>0.$ Then the $D$ matrix reduces to
$$D=\left(
\begin{array}{ccc}
 d_{11}&\frac{\alpha_{12}}{b_{12}}&d_{13}\\
  d_{21}&\frac{\alpha_{22}}{b_{22}}&d_{23}\\
    \end{array}
\right).$$
 Again, from the compatibility condition, we know that if $A$ and $B$  are compatible, then rank($D$)$>$ 1. However, in this case $\text{rank(D)}\leq \text{min(2,3)}=2.$ So, in this case we must have rank ($D$)=2. Thus, any $(2\times 2)$ determinant will be non-vanishing  (for all admissible choices of
 $0<(\alpha_{12},\alpha_{22})\leq 1$ with $ \alpha_{12}+\alpha_{22}=1$ ) means that when
 $$ \frac{\alpha_{22}}{b_{22}}d_{11}-d_{21}\frac{\alpha_{12}}{b_{12}}\neq 0$$ and

 $$ \frac{\alpha_{12}}{b_{12}}d_{23}-d_{13}\frac{\alpha_{22}}{b_{22}}\neq 0,$$ the two matrices $A$ and $B$ will be incompatible.

 Again, let us consider the  compatibility set-up corresponding to which we have the  equation
 \begin{equation*}
 a_{ij}\sum_{s=1}^{I}b_{sj}\eta_{s}-b_{ij}\eta_{i}=0,
 \end{equation*}

 and  $\eta_{i}$ are the marginals corresponding to the variable $X.$ The above equation can be written in terms of a system of  equations as mentioned earlier, for which the $D$ matrix reduces to
\begin{center}
$\begin{bmatrix}
 b_{11}(a_{11}-1)& a_{11}b_{21}&a_{11}b_{31}\\
  b_{12}(\alpha_{12}-1)&\alpha_{12}b_{22}& \alpha_{12}b_{32}\\
   a_{21}b_{11} & b_{21}(a_{21}-1) &a_{21}b_{31} \\
   \alpha_{22}b_{12}& b_{22}(\alpha_{22}-1)&\alpha_{22}b_{32}\\
     a_{31}b_{11} &a_{31}b_{21} &b_{31}(a_{31}-1)\\
    a_{32}b_{12} &a_{32}b_{22} & b_{32}(a_{32}-1).\\
    \end{bmatrix}$
\end{center}

So,  if $A$ and $B$ are not compatible and we instead consider the concept of $\varepsilon$-compatibility,  then according to Arnold et al. (1999), the system of inequalities can be written as
\begin{equation}
\begin{bmatrix}
 b_{11}(a_{11}-1)& a_{11}b_{21}&a_{11}b_{31}\\
  b_{12}(\alpha_{12}-1)&\alpha_{12}b_{22}& \alpha_{12}b_{32}\\
   a_{21}b_{11} & b_{21}(a_{21}-1) &a_{21}b_{31} \\
   \alpha_{22}b_{12}& b_{22}(\alpha_{22}-1)&\alpha_{22}b_{32}\\
     a_{31}b_{11} &a_{31}b_{21} &b_{31}(a_{31}-1)\\
    a_{32}b_{12} &a_{32}b_{22} & b_{32}(a_{32}-1)\\
    \end{bmatrix}  \underline{\eta}
    \leq \underline{\varepsilon}.
\end{equation}

There will be another six equations which will be  exactly the same,  but with only difference in their sign. Note that in this situation we have a system of equations with $3$ unknowns (including those  constraints on $\eta$ and $\alpha_{ij}).$ Let us try to find out those values by considering equality in the previous set of equations.  Thus, from  (37), we get

\begin{equation}
 a_{31}b_{11}\eta_{1} +a_{31}b_{21}\eta_{2} +b_{31}(a_{31}-1)\eta_{3}=\varepsilon,
 \end{equation}

\begin{equation}
 a_{32}b_{12}\eta_{1} +a_{32}b_{22}\eta_{2} +b_{32}(a_{32}-1)\eta_{3}=\varepsilon,
 \end{equation}

and
\begin{equation}
  \eta_{1}+\eta_{2}+\eta_{3}=1.
\end{equation}

 So, after some algebra, we get (by substituting the expression for $\eta_{3}=1-(\eta_{1}+\eta_{2}),$ in (38) and (39), and solving)
 \begin{equation}
 \eta_{2}=\frac{\varepsilon d_{12}+d_{22}}{d_{11}},
 \end{equation}

 where
 $$d_{11}=(a_{31}b_{21}-b_{31}(a_{31}-1))(a_{32}b_{12}-b_{32}(a_{32}-1))-(a_{32}b_{22}-b_{32}(a_{32}-1))(a_{31}b_{11}-b_{31}(a_{31}-1)),$$

 $$d_{12}=a_{32}b_{12}-b_{32}(a_{32}-1)-a_{31}b_{11}+b_{31}(a_{31}-1),$$
 and
 $$ d_{22}=b_{32}(a_{32}-1)(a_{31}b_{11}-b_{31}(a_{31}-1))-b_{31}(a_{31}-1)(a_{32}b_{12})(a_{32}b_{12}-b_{32}(a_{32}-1)).$$

\noindent Furthermore,

 \begin{equation}
 \eta_{1}=\frac{d_{11}[\varepsilon-(a_{31}-1)b_{31}]-[a_{31}b_{21}-b_{31}(a_{31}-1)](\varepsilon d_{12}+d_{22})}
 {a_{31}b_{11}-b_{31}(a_{31}-1)}.
 \end{equation}

  \noindent The estimated value of  $\eta_{3}$ can be found by substituting the values from (41) and (42) into (40). Thus,  estimates of the unknown values of $\alpha_{12}$ and  $\alpha_{22}$ will be

  \begin{equation*}
  \alpha_{12}=\frac{\varepsilon+b_{12}\eta_{1}}{b_{12}\eta_{1}+b_{21}\eta_{2}+b_{31}\eta_{3}}
  \end{equation*}

  and $$\alpha_{22}=1-\alpha_{12}.$$

  \bigskip

 \end{enumerate}

\subsection{Compatibility in the General Case}
We now will discuss the problem of  compatibility when the dimension of the two matrices $A$ and $B$ is of the order ($I\times J$). In this case, we restrict to  the situation when there are two elements unknown only in $A,$ while in matrix $B$ all the elements are known. Let us consider  in matrix $A$ in the $l_{1}-th$ column ($1\leq l_{1} \leq J$),  two elements are unknown and  they appear at  $i_{1}-th$  and $i_{2}-th$  rows, where $(1\leq (i_{1},i_{2})\leq I).$  Since column sums of $A$ add up to 1, we can write, considering the unknown elements to be denoted by $\alpha_{ij}'$s,

\begin{equation}
\alpha_{i_{1}l_{1}}+\alpha_{i_{2}l_{1}}+\sum_{k\neq i_{1},i_{2}}a_{kj}=1, \forall (j,l_{1},l_{1})=1(1)J,
\end{equation}
while for $B$ all the rows  add up to 1,  all the elements are known, and that
$$ \sum_{j=1}^{J}b_{ij}=1, \forall i=1(1)I.$$
Now if we have the information that the matrices $A$ and $B$ are compatible, then  we can write
\begin{gather}
\alpha_{i_{1}l_{1}}\Bigg[\sum_{s=1}^{I}b_{sl_{1}}\eta_{s}\Bigg]-b_{i_{1}l_{1}}\eta_{i_{1}}=0,\\
\Rightarrow \alpha_{i_{1}l_{1}}\Bigg[b_{i_{1}l_{1}}\eta_{i_{1}}+\sum_{\substack{s=1\\s\neq i_{1}}}b_{sl_{1}} \eta_{s}\Bigg]
-b_{i_{1}l_{1}}\eta_{i_{1}}=0.
\end{gather}
Hence, the unknown value of $\alpha_{i_{1}l_{1}}$ will be given by
$$\alpha_{i_{1}l_{1}}=\frac{b_{i_{1}l_{1}}\eta_{i_{1}}}{b_{i_{1}l_{1}}\eta_{i_{1}}+\sum_{\substack{s=1\\s\neq i_{1}}}b_{sl_{1}} \eta_{s}}.$$

\noindent We can  then  write
$$ \alpha_{i_{2}l_{1}}=1-\sum_{k\neq i_{1},i_{2}}a_{kj}-\alpha_{i_{1}l_{1}}
= 1-\sum_{k\neq i_{1},i_{2}}a_{kj}-\frac{b_{i_{1}l_{1}}\eta_{i_{1}}}{b_{i_{1}l_{1}}\eta_{i_{1}}+\sum_{\substack{s=1\\s\neq i_{1}}}b_{sl_{1}} \eta_{s}}.$$
However, the solution for  $\eta_{i}'$ s can be obtained by using any set of $(I-1)$ equations, since
$\sum_{i=1}^{I}\eta_{i}=1.$

\section{Concluding Remarks}
The  search for a compatible $P$ in terms of  equations subject to inequality constraints is based on the fact that we really need to find one compatible marginal, say,  corresponding to the random variable $X$, and we consider the fact that when this is combined with $B$ will give us $P.$ Compatible conditional and marginal specifications of distributions are of fundamental importance in modeling. Moreover, in Bayesian prior elicitation context, inconsistent conditional specifications are to be expected. In such situations, interest will center on most nearly compatible distributions (see Arnold et al. (1999), and Ghosh and Balakrishnan (2013)).
In the finite discrete case, a variety of compatibility conditions can be derived. In this article, we have discussed in detail the problem of compatibility in  the context  mentioned earlier by identifying it as a  programming problem and have developed a rank--based criterion.  Although we have shown that the rank of the matrix  (whose elements are constructed from the two given matrices $A$ and $B$) under compatibility will be $I-1$, for a ($3\times 3$), and ($4\times 4$) case, the result is true for any dimension of the two given matrices. A significant amount work here  draws heavily from Arnold et al. (1999) and Arnold and Gokhale (1998). Also, we have provided a discussion  on the topic of compatibility when we have missing elements in either $A$ or  $B,$ or in both. It has been observed that for a given $A$ and $B,$ under compatibility, the choices for the missing elements are unique.  In addition,  we have discussed in this context what would be the possible choices of those missing elements when we have the information that the two matrices are incompatible. However, for a general case when the dimension of the matrix $D$ is $\left(IJ \times I\right),$ the strategy  discussed here
will be quite challenging in identifying the solution for the unknown
elements either in any of the conditional probability matrices A and B
or in both. Also, when we have elements missing in $A$ and $B$
in different positions, then the procedure will result in solving a set
of $IJ$ number of equations which is cumbersome and quite difficult to
handle. In such a situation, one may consider the concept of
compatibility under rank one criterion as proposed by Arnold et al. (2001).
One interesting question that may arise here is how can we extend the above technique under compatibility when there exists more than two conditional matrices in the discrete case,  i.e.,  if we are given three matrices $A,$ and $B$ and $C,$ where\\

 \noindent $A$ is the conditional probability matrix of (say) $X,$ given $Y$ and $Z,$\\
 $B$ is the conditional probability matrix of (say) $Y,$ given $X$ and $Z,$\\
 $C$ is the conditional probability matrix of (say) $Z,$ given $X$ and $Y.$\\

\noindent Furthermore, what would happen in the situation (under compatibility) when some of the elements  are unknown in  any of $A,$  $B$ or $C,$ or in all of them? Such questions require a careful  study of the concept of  compatibility. We are currently looking into these issues and hope to report the findings in a  future paper.

\end{document}